\documentclass[11pt]{article}
\usepackage[colorlinks=true,
linkcolor=webgreen,
filecolor=webbrown,
citecolor=webgreen]{hyperref}

\usepackage{amssymb,psfig,epsfig}

\definecolor{webgreen}{rgb}{0,.5,0}
\definecolor{webbrown}{rgb}{.6,0,0}

\newcommand{\hsp}{\hspace*{\parindent}}
\newcommand{\eqn}[1]{(\ref{#1})}
\newcommand{\eeq}{\end{equation}}
\newcommand{\beql}[1]{\begin{equation}\label{#1}}

\newcommand{\la}{\lambda}

\newcommand{\sV}{{\mathcal V}}
\newcommand{\RR}{{\mathbb R}}

\newcommand{\ZZ}{{\mathbb Z}}
\newcommand{\sP}{{\mathcal P}}

\makeatletter
\def\@sect#1#2#3#4#5#6[#7]#8{\ifnum #2>\c@secnumdepth
     \def\@svsec{}\else
     \refstepcounter{#1}\edef\@svsec{\csname the#1\endcsname.\hskip .75em }\fi
     \@tempskipa #5\relax
      \ifdim \@tempskipa>\z@
        \begingroup #6\relax
          \@hangfrom{\hskip #3\relax\@svsec}{\interlinepenalty \@M #8\par}%
        \endgroup
       \csname #1mark\endcsname{#7}\addcontentsline
         {toc}{#1}{\ifnum #2>\c@secnumdepth \else
                      \protect\numberline{\csname the#1\endcsname}\fi
                    #7}\else
        \def\@svsechd{#6\hskip #3\@svsec #8\csname #1mark\endcsname
                      {#7}\addcontentsline
                           {toc}{#1}{\ifnum #2>\c@secnumdepth \else
                             \protect\numberline{\csname the#1\endcsname}\fi
                       #7}}\fi
     \@xsect{#5}}
\def\@begintheorem#1#2{\it \trivlist \item[\hskip \labelsep{\bf #1\ #2.}]}

\def\section{\@startsection {section}{1}{\z@}{-3.5ex plus -1ex minus 
 -.2ex}{2.3ex plus .2ex}{\normalsize\bf}}
\makeatother
\makeatletter
\def\subsection{\@startsection {subsection}{1}{\z@}{-3.5ex plus -1ex minus
 -.2ex}{2.3ex plus .2ex}{\normalsize\bf}}

\makeatother

\begin{document}
\begin{center}
{\large {\bf A Zador-Like Formula for Quantizers Based on Periodic Tilings}} \\
\vspace{1.5\baselineskip}
N. J. A. Sloane and Vinay A. Vaishampayan \\ [+.2in]
Information Sciences Research Center \\
AT\&T Shannon Lab \\
Florham Park, New Jersey 07932-0971, USA \\
Email: \href{mailto:njas@research.att.com}{njas@research.att.com},
\href{mailto:vinay@research.att.com}{vinay@research.att.com} \\
\vspace{1.5\baselineskip}
July 16, 2002 \\
\vspace{1.5\baselineskip}

{\bf ABSTRACT}
\vspace{.5\baselineskip}
\end{center}
\setlength{\baselineskip}{1.5\baselineskip}

We consider Zador's asymptotic formula for the distortion-rate function for a variable-rate vector quantizer in the high-rate case.
This formula involves the differential entropy of the source, the rate of the quantizer in bits per sample, and a coefficient $G$ which depends on the geometry of the quantizer but is independent of the source.
We give an explicit formula for $G$ in the case when the quantizing regions
form a periodic tiling of $n$-dimensional space, in terms of the volumes and second moments of the
Voronoi cells.
As an application we show, extending earlier work of Kashyap and Neuhoff, that even a variable-rate
three-dimensional quantizer based on the ``A15'' structure is still inferior to a quantizer based on the body-centered cubic lattice.
We also determine the smallest covering radius of such a structure.



\section{Introduction}
\hsp
Zador gave two asymptotic formulae for the distortion-rate function for a vector quantizer in $\RR^n$, depending on whether a fixed- or variable-rate code is
used to transmit which cell the sample point belongs to.
These are treated as cases (A) and (B) in \cite{Z2}.
See also Equations (19) and (20) in \cite{GrNe98}.
In the present note we are concerned with the variable-rate case.
Zador's formula for this case states that
the distortion-rate function
has the form
\beql{Eq1}
\delta (R) \cong G 2^{2h(X)} 2^{-2R} ~,
\eeq
where $\delta (R)$ is the average squared error per symbol,
$h(X)$ is the differential entropy per dimension of the source $X$,
$R$ bits/symbol is the rate of the quantizer,
and $G$ depends on the quantizer but not on the source.
This is from page 4 of \cite{Z2} and Equation (20) in \cite{GrNe98}.
Note that, as pointed out in \cite{GrNe98}, the version of Eq. (1)
given in \cite{Z3} is incorrect.
For further information see also \cite{Z1} and \cite{Gray}, Chapter 5.
The value of $G$ can therefore be used to compare different quantizers.
To calculate $G$ we may choose any convenient distribution for the source $X$.

Suppose first that the quantizer points form a lattice $\Lambda \subseteq \RR^n$,
with Voronoi cell $\sV$ of volume $V = \sqrt{\det \Lambda}$.
To calculate $G$ we assume that $X$ is uniformly distributed over $\sV$.
Then $G$ becomes
\beql{Eq2}
G = \frac{1}{n} \frac{\int_\sV \| x \|^2 dx}{V^{1+2/n}} ~,
\eeq
the familiar expression for the average mean squared error per dimension of a lattice quantizer
(\cite{CSV}, \cite[Chap. 2]{SPLAG}, \cite[Chap. 10]{GG92}, \cite[Chap. 5]{Gray}).
In this case all the quantizing cells have the same volume, and the fixed-rate and variable-rate versions of Zador's formula coincide.

The main purpose of this note is to put on record an analogue of \eqn{Eq2} (see Eq. \eqn{EqU5} and \eqn{EqU5B}) for the variable-rate Zador formula which applies to the case
when the quantizing regions form a periodic tiling of $\RR^n$, such as that formed by the Voronoi cells for a union of a finite number of translates of
a lattice.
This question arose in a recent investigation of quantizers that are based on writing a lattice as an intersection of several simpler lattices \cite{SlBe}.

A second motivation for our work was a recent paper of Kashyap and Neuhoff \cite{KaNe}, which makes use of a fixed-rate analogue of \eqn{Eq2} for periodic lattice quantizers (see Eq. \eqn{EqK1} and \eqn{EqKB} below).
One of the goals of \cite{KaNe} was to see if the A15 arrangement of points
in $\RR^3$  that has recently arisen in several different contexts
(\cite{KHPC},  \cite{LaSh}, \cite{W}, \cite{WP94}) could produce a
better three-dimensional quantizer than the body-centered cubic (bcc) lattice.
The latter is known to be the best {\em lattice} quantizer in $\RR^3$ \cite{BS1}, but the question of the existence of a better nonlattice quantizer remains open.

Kashyap and Neuhoff \cite{KaNe} found that with their figure of merit the best
quantizer based on the A15 arrangement is inferior to the bcc lattice.
In Section 3 we repeat the comparison using our figure of merit.
Now a different version of the A15 quantizer is best, but
is still inferior to the bcc lattice.

Other applications of our formula will be found in \cite{SlBe}.

Another open problem in three-dimensional geometry is to determine the best covering of $\RR^3$ by equal (overlapping) spheres.
In Section 3 we determine the smallest covering radius that can be achieved with the A15 structure.
This is also
(slightly) worse than
that of the bcc lattice.

\section{Periodic quantizers}
\hsp
Consider a vector quantizer in $\RR^n$ in which the quantizing regions form a periodic tiling.
Let $\sV$ be a minimal periodic unit or tile for the tiling, and let $\sP_1, \ldots, \sP_k$ be a list of the different polytopes occurring among the quantizing regions.
Let $c_i$ be the centroid of $\sP_i$.
In order to determine $G$ we assume
the source $X$ is uniformly distributed over $\sV$, and let $p_i$ $(i=1, \ldots, k)$ be the probability that $X$ is in a cell of type $\sP_i$.
Also let $N_i = p_i V /V_i$, where $V_i = \mbox{vol}\, \sP_i$ and $V= \mbox{vol} \, \sV$.
Then $N_i$ is the number of cells of type $\sP_i$ per copy of $\sV$, and
\beql{EqU0}
V = N_1 V_1 + \cdots + N_k V_k ~.
\eeq

For example, let $\Lambda \subseteq \RR^n$ be a lattice and let $a_1 + \Lambda, \ldots, a_r + \Lambda$ $(a_i \in \RR^n )$ be distinct translates.
Then the Voronoi cells for the union of the points $a_i + \Lambda$ $(i=1, \ldots, r )$ form a quantizer of the type considered here.

We now apply \eqn{Eq1}.
The left-hand side is the normalized mean squared error per dimension $U/(n V)$, where
\beql{EqU1}
U = \sum_{i=1}^k N_i U_i = V \sum_{i=1}^k \frac{p_i U_i}{V_i}
\eeq
and
\beql{EqU1a}
U_i = \int_{\sP_i} \| x-c_i \|^2 dx
\eeq
is the unnormalized mean squared error over a cell of type $\sP_i$.

The differential entropy per dimension is
\beql{EqU3}
h(X) = \frac{1}{n} \log_2 V , ~~\mbox{so}~~
2^{2h (X)} = V^{2/n} ~.
\eeq
It remains to calculate the rate $R$ of the quantizer.
Observe that we need
$H(p_1, \ldots, p_k) = - \sum_{i=1}^k p_i \log_2 p_i$ bits to specify the type of cell to which the quantized point belongs, and a further
$\sum_{i=1}^k p_i \log_2 N_i = \sum_{i=1}^k p_i \log_2 (p_i V/V_i )$ bits to specify the particular one of the $N_i$ cells of that type.
This requires a total of $\log_2 V - \sum_{i=1}^k p_i \log_2 V_i$ bits, and
then $R$ is this quantity divided by $n$, so that
\beql{EqU4}
2^{-2R} = V^{-2/n} \prod_{i=1}^k V_i^{2p_i /n} ~.
\eeq

We substitute \eqn{EqU1}, \eqn{EqU3}, \eqn{EqU4} into \eqn{Eq1} to get our expression
\beql{EqU5}
G = \frac{\sum_{i=1}^k \frac{p_i U_i}{V_i}}{n \left( \prod_{i=1}^k V_i^{p_i}\right)^{\frac{2}{n}}} ~,
\eeq
for the average mean squared error per dimension
using variable-rate coding.
The numerator of \eqn{EqU5} is equal to $U/V$ (see \eqn{EqU1}),
so we may rewrite \eqn{EqU5} as
\beql{EqU5B}
G = \frac{U}{n V \prod_{i=1}^k V_i^{2p_i /n}}
\eeq

In contrast, the expression given by Kashyap and Neuhoff \cite{KaNe} for fixed-rate coding is the following.
Suppose the basic tile $\sV$ contains $L$ cells, $Q_1, \ldots, Q_L$, not assumed to be
distinct, where $Q_i$ has volume $V_i$ and unnormalized
second moment $U_i$ (as in \eqn{EqU1a}).
Then their expression is
\beql{EqK1}
G = \frac{L^{2/n} \sum_{i=1}^L U_i}{nV^{1+2/n}} ~,
\eeq
where $V = \mbox{vol} \, \sV$.
The following justification is equivalent to the one in \cite{KaNe}, but
clarifies the difference between our two approaches.
Consider a large region of space, $B$, which is partitioned into $\lambda$ copies of $\sV$, and let $X$ be uniformly distributed over $B$.
The left-hand side of \eqn{Eq1} is
\beql{EqK2}
\frac{\lambda \sum_{i=1}^L U_i}{n \, \mbox{vol} \, B} = \frac{\sum_{i=1}^L U_i}{n V} ~.
\eeq
The differential entropy per dimension is
\beql{EqK3}
h(X) = \frac{1}{n} \log_2 ( \mbox{vol} \, B ), ~~\mbox{so}~~
2^{2h (X)} = (\la V )^{2/n} ~.
\eeq
To compute the rate, it may be argued that $\log_2 \la$ bits are required to specify the copy of $\sV$, and $\log_2 L$ bits to specify which of the $Q_i$'s the point belongs to.
This requires a total of $\log_2 \la L$ bits, so $R = (1/n) \log_2 \la L$, and
\beql{EqK4}
2^{-2R} = (\la L)^{-2/n} ~.
\eeq

Arguing as before, we substitute \eqn{EqK2}--\eqn{EqK4} into \eqn{Eq1},
which gives \eqn{EqK1}.

If there are $k$ different types of cells among the $Q_i$, with the $j$-th cell occurring $N_j$ times, then \eqn{EqK1} can be written as
\beql{EqKB}
G = \frac{\left( \sum_{i=1}^k N_i \right)^{2/n} U}{nV^{1+2/n}} ~.
\eeq

The ratio of the two expressions, \eqn{EqU5B} divided by \eqn{EqKB},
can be written as
\beql{EqRatio}
2^{ -\frac{2}{n} \left( \log_2 L ~-~ \sum_{i=i}^k p_i \log_2 (p_i/N_i)  \right) } ~.
\eeq
Both formulae, \eqn{EqU5B} and \eqn{EqKB}, depend only on the geometry of the quantizer.
The difference between the two expressions arises because, as long as the cells do not all have the same volume, variable-length coding can take advantage of the different cell probabilities to reduce the overall rate.

If all cells have the same volume then \eqn{EqU5B} and \eqn{EqKB} coincide, and if there is only type of cell then they both reduce to \eqn{Eq2}.

In general, the fact that \eqn{EqU5B} is less than equal to \eqn{EqKB} can be shown directly.
After canceling some common factors and rearranging, we must show that
$$
\prod_{j=1}^k V_j^{p_j} \ge \frac{1}{\sum_{i=1}^k \frac{p_i}{V_i}} ~,
$$
and this follows from the geometric-mean harmonic-mean inequality.

\section{Quantizers based on the A15 structure}
\hsp
The A15 structure was discovered by Kasper et~al. \cite{KHPC} in the clathrate compound ${\rm Na_8 Si_{46}}$.
It was used by Weaire and Phelan (\cite{WP94}; see also \cite{W}) to construct a counter-example to
Kelvin's conjecture on minimal surface soap films.
It was also used by Lagarias and Shor \cite{LaSh} to construct counter-examples
to Keller's conjecture in dimensions 10 and above.
Since the soap film problem attempts to find a partition of space
into cells which are good approximations to spheres,
it was therefore natural to ask if the A15 structure could also lead to a record-breaking
quantizer in three dimensions.

The A15 structure may be defined as the union of eight translates of the cubic lattice $4 \ZZ^3$ by the vectors $(0,0,0)$, $(2,2,2)$, $(0 \pm 1, 2)$, $(2, 0, \pm 1 )$, $(\pm 1, 2, 0)$.
(The first two translates alone give the bcc lattice.)

There are two types of points in this structure, the {\em even} points in which all coordinates are even, and the {\em odd} points in which some coordinate is odd.
There are isometries of $\RR^3$ which map A15 to itself and
act transitively on the even points and on the odd points. However, 
even points are not equivalent to odd points.

Weaire and Phelan consider a weighted Voronoi decomposition of $\RR^3$:
walls between points of the same type occur along perpendicular bisectors,
but a wall between an even point $E$ and an odd point $D$ is defined by the plane
$$(X-E) \cdot (D-E) = \mu (D-E) \cdot (D-E) ~,
$$
where $\mu$ is a weighting factor to be determined.
For $0 < \mu < 3/5$ there are two types of Voronoi cells:
12-sided polyhedra $\sP_1$ centered at the even points and 14-sided polyhedra
$\sP_2$ centered at the odd points.
Weaire and Phelan use the ``Surface Evolver'' computer program of Brakke \cite{Bra} to perturb the polyhedra so that they have equal volumes and minimal total surface area.

Kashyap and Neuhoff \cite{KaNe} use the weighted Voronoi decomposition based
on the A15 structure, but
adjust $\mu$ to give the smallest value of \eqn{EqK1}.
The formula are simpler if we set $\mu = 2 \alpha /5$.
Then Kashyap and Neuhoff find (and we have confirmed) that
\begin{eqnarray*}
V_1 & = & \mbox{vol}\, \sP_1 = 4\alpha^3 ~, \\
V_2 & = & \mbox{vol} \, \sP_2 = \frac{4}{3} (8 - \alpha^3 ) ~, \\
U_1 & = & \frac{71}{30} \alpha^5 ~, \\
U_2 & = & \frac{1}{90}
(1200 - 600 \alpha^3 + 360 \alpha^4 - 71 \alpha^5 ) \,.
\end{eqnarray*}
The fundamental tile $\sV$ has volume $V=64$, and $L=8$, $N_1=2$, $N_2=6$, $p_1 = \alpha^3 /8$, $p_2 = 1- \alpha^3/8$.
Then \eqn{EqK1} becomes
\beql{EqL1}
\frac{8^{2/3} (2U_1 + 6 U_2 )}{3 \cdot 64^{5/3}} = \frac{1}{96}
(10 - 5 \alpha^3 + 3\alpha^4 ) \,.
\eeq
This has a minimal value of 0.07873535... at $\alpha = 5/4$.
For comparison, the values of $G$ for the bcc and fcc (face-centered cubic)
lattices are $19/(192 \cdot 2^{1/3} ) = 0.07854328\ldots$
and
$2^{-11/3} = 0.07874507\ldots$, respectively (\cite{BS1}, \cite{CSV}, \cite{SPLAG}).

Using our formula \eqn{EqU5} we find that
\beql{EqL2}
G = \frac{2U_1 + 6 U_2}{192 (V_1^{\alpha^3/8} V_2^{1- \alpha^3/8})^{2/3}}
\eeq
which has a minimal value of 0.07872741$\ldots$ at $\alpha = 1.2401 \ldots$.
This is slightly better, but still inferior to the bcc lattice.

Another unsolved problem is to find the best covering
of $\RR ^3$ by overlapping spheres (cf. \cite[Chap. 2, Table 2.1]{SPLAG}).
The bcc lattice is the best {\em lattice} covering, with thickness $1.4635 \ldots$,
but the question of the existence of a better nonlattice covering remains open.
We find that the covering radius of the weighted Voronoi decomposition of A15
is minimized at $\alpha = 5/4$, which gives a thickness of
$$\frac{125 \sqrt(3) \pi}{432} = 1.5745  \dots ~,$$
just slightly worse than that of the bcc lattice
(but again better than the fcc lattice, which has thickness
$2.0944\ldots$).

\subsection*{Acknowledgment}
We thank the referees for some helpful comments.

\end{document}